\documentclass[12pt]{article}
\usepackage[final]{epsfig}
\usepackage{graphics}
\usepackage{amsmath}
\usepackage{amsfonts}
\usepackage{latexsym}
\usepackage{amssymb}
\usepackage{amsthm}
\usepackage{graphicx}
\usepackage{url}
\usepackage{epstopdf}
\usepackage[T2A,T1]{fontenc}
\newcommand{\YA}{\mbox{\usefont{T2A}{\rmdefault}{m}{n}\CYRYA}}
\newcommand{\FF}{\mbox{\usefont{T2A}{\rmdefault}{m}{n}\CYRF}}
\newcommand{\aaa}{\mbox{\usefont{T2A}{\rmdefault}{m}{n}\cyra}}

\newtheorem{lemma}{Lemma}[section]

\newtheorem{theorem}{Theorem}
\newtheorem{corollary}[lemma]{Corollary}
\newtheorem{conjecture}{Conjecture}

\theoremstyle{definition}
\newtheorem{remark}[lemma]{Remark}
\newtheorem{definition}{Definition}

\newcommand{\proofend}{$\Box$\bigskip}

\newcommand{\R}{{\mathbb R}}

\def\proof{\paragraph{Proof.}}
\newcommand{\g}{\gamma}
\newcommand{\G}{\Gamma}
\def\da{\partial a}
\def\wt#1{\widetilde#1}
\def\db{\partial b}
\def\dx{\partial x}
\def\dy{\partial y}
\def\mca{\mathcal A}
\def\la{\lambda}
\def\rr{\mathbb R}
\def\mcd{\mathcal D}
\def\const{\mathrm{const}}

\title{Four equivalent properties of integrable billiards}

\author{Alexei Glutsyuk\footnote{
CNRS, France (UMR 5669 (UMPA, ENS de Lyon) and UMI 2615 (Interdisciplinary Scientific Center J.-V.Poncelet))
\newline\indent
National Research University Higher School of Economics (HSE), Moscow, Russia
\newline\indent
E-mail: 
aglutsyu@ens-lyon.fr}
\and
Ivan Izmestiev\footnote{
TU Wien,
Wiedner Hauptstra{\ss}e 8-10/104,
1040 Vienna, Austria \newline\indent
E-mail:
izmestiev@dmg.tuwien.ac.at}
\and
Serge Tabachnikov\footnote{
Department of Mathematics,
Penn State University,
University Park, PA 16802, USA
\newline\indent
E-mail: 
tabachni@math.psu.edu}
}

\date{}

\begin{document}

\maketitle

\vspace{-.3cm}
\begin{abstract}
By a classical result of Darboux, a foliation of a Riemannian surface has the Graves property
(also known as the strong evolution property)
if and only if the foliation comes from a Liouville net. A similar result of Blaschke says that a pair of orthogonal foliations has the Ivory property if and only if they form a Liouville net.

Let us say that a geodesically convex curve on a Riemannian surface has the Poritsky property if it can be parametrized in such a way that all of its string diffeomorphisms are shifts with respect to this parameter. In 1950, Poritsky has shown that the only closed plane curves with this property are ellipses.

In the present article we show that a curve on a Riemannian surface has the Poritsky property if and only if it is a coordinate curve of a Liouville net.
We also recall Blaschke's derivation of the Liouville property from the Ivory property and his proof of Weihnacht's theorem: the only Liouville nets in the plane are nets of confocal conics and their degenerations.

This suggests the following generalization of Birkhoff's conjecture: If an interior neighborhood of a closed geodesically convex curve on a Riemannian surface is foliated by billiard caustics, then the metric in the neighborhood is Liouville, and the curve is one of the coordinate lines.
\end{abstract}

\section{Introduction} \label{sect:intro}

\subsection{Billiards, caustics, and string construction}
The billiard dynamical system describes the motion of a free particle (billiard ball) inside a domain (billiard table): the particle moves with a constant velocity and reflects elastically from the boundary so that the angle of incidence equals angle of reflection. This law is familiar from geometrical optics. 

We refer to the books \cite{CM06,KT91,Tab95,Tab05} for the theory of mathematical billiards and, in particular, to various properties that we discuss below. If the billiard table is a convex body, one 
can reduce the continuous time billiard system to a map that acts on oriented lines that intersect the billiard table (one may think of them as rays of light): this {\it billiard ball map} takes the incoming trajectory of the billiard ball to the outgoing one. In this paper, we will study billiard tables with a strictly convex smooth boundary. 

A billiard {\it caustic} is a curve with the property that if a segment of billiard trajectory is tangent to it, then the reflected segment is also tangent to it. That is, a caustic corresponds to an invariant curve of the billiard ball map. Although caustics may have singularities (generically, semi-cubical cusps), we will consider only smooth closed convex caustics.

It is a deep result of the KAM theory that if the boundary of the billiard table is  strictly  convex and sufficiently smooth then, arbitrarily close to the boundary, there exist caustics; furthermore, these caustics occupy a set of positive measure \cite{Laz}. In general, these caustics do not form a foliation, and there are slits between them.  

\begin{figure}[hbtp]
\centering
\includegraphics[height=2in]{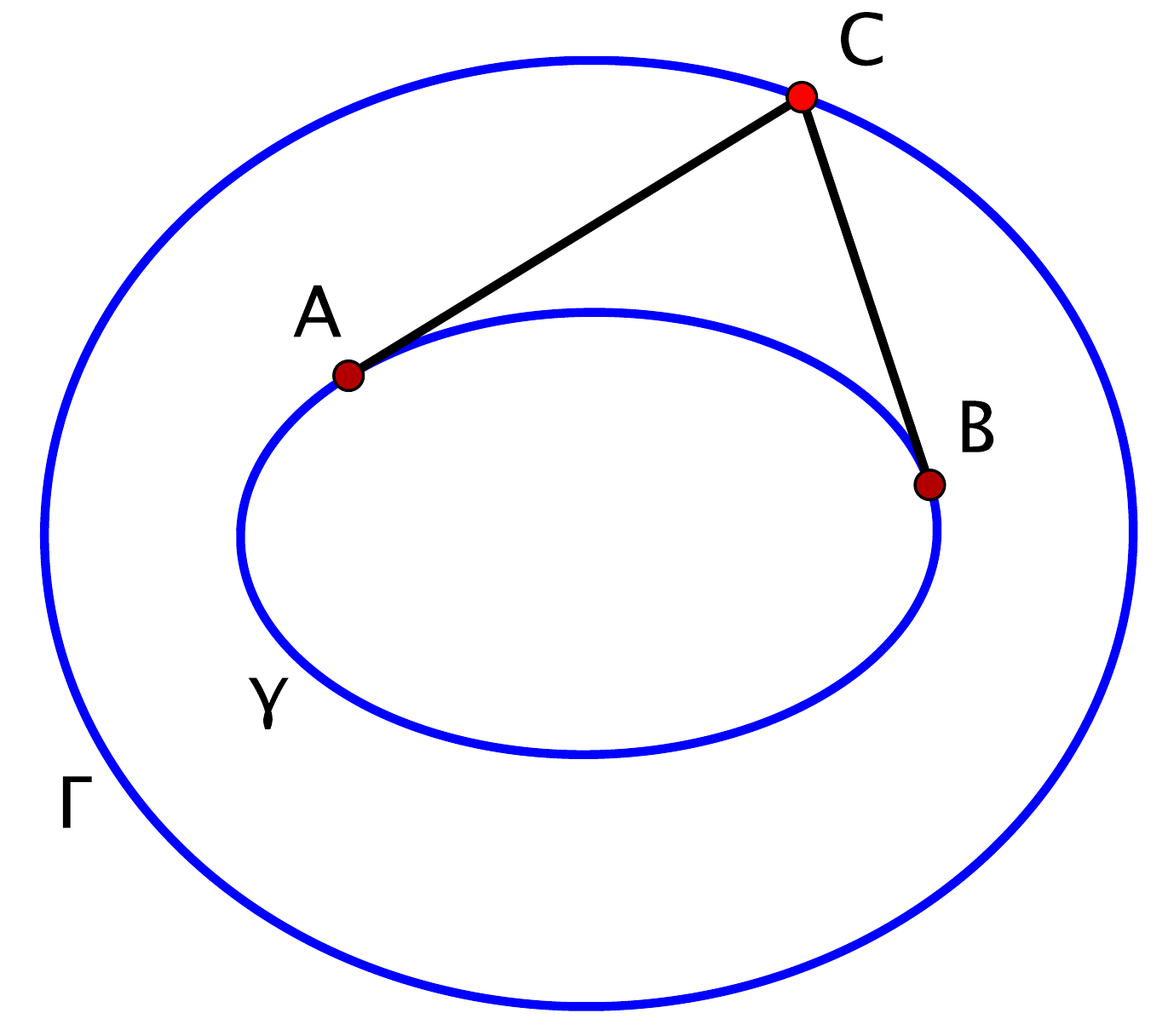}
\caption{The string construction on $\g$ yields the billiard curve $\G$.}
\label{string}
\end{figure}

Given a billiard table, there is no easy way to find a caustic, but the converse problem, to construct a table from a caustic, is solved by the following {\it string construction}. Let $\g$ be a caustic. Wrap a non-stretchable string around $\g$, pull it tight at a point, and move the point around keeping the string tight. The trajectory of the point is a curve $\G$, the boundary of a billiard table for which $\g$ serves as a caustic, see Figure \ref{string}. This construction yields a 1-parameter family of curves $\G$, the parameter being the length of the string.

Once the length of the string is chosen, the string construction defines a diffeomorphism of the curve $\gamma$ that takes one tangency point of the string with the curve to the other point. Assuming that the curve $\gamma$ in Figure \ref{string} is oriented clockwise, this diffeomorphism takes point $A$ to $B$. Let us call these maps {\it string diffeomorphisms}.
Varying the length of the string, one has a 1-parameter family of string diffeomorphisms of the same curve.

\subsection{Confocal conics}
Beyond the nearly trivial case of a circle (whose caustics are concentric circles), the most classical example is the billiard inside an ellipse. Let us recall the salient features of this system; see, e.g., the above mentioned books and \cite{Iz17}.

The billiard inside an ellipse is completely integrable. Specifically, 
the caustics of the billiard inside an ellipse are confocal ellipses: if a billiard trajectory does not intersect the segment between the foci, it is tangent to a unique confocal ellipse, a caustic. A billiard trajectory that passes through a focus reflects to the other focus. 
If a billiard trajectory intersects the segment between the foci, it is tangent to a unique confocal hyperbola, which is also a caustic. 
 
\begin{figure}[hbtp]
\centering
\includegraphics[height=1.7in]{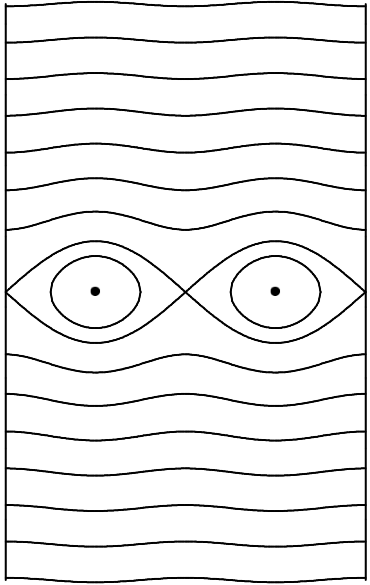}
\caption{Phase portrait of the billiard ball map in an ellipse.}
\label{phase}
\end{figure}

The difference between two types of caustics is that the ellipses correspond to non-contractible invariant circles of the billiard ball map, whereas the hyperbolas correspond to contractible invariant curves. The phase space of the billiard is a cylinder, see Figure \ref{phase}. 

As a consequence, one has the Graves theorem: {\it wrapping a closed non-stretchable string around an ellipse produces a confocal ellipse.} Equivalently, confocal ellipses satisfy what is known as the strong evolution property: a caustic of a caustic is again a caustic. 

This prompts the following definition.

\begin{definition} \label{def:Graves}
{\rm 
An annulus foliated by simple closed convex curves has the {\it Graves property} if, for every pair of nested leaves, the outer one is obtained by the string construction from the inner one. In particular, this applies to the boundary curves. 
}
\end{definition}

\begin{figure}[hbtp]
\centering
\includegraphics[height=1.8in]{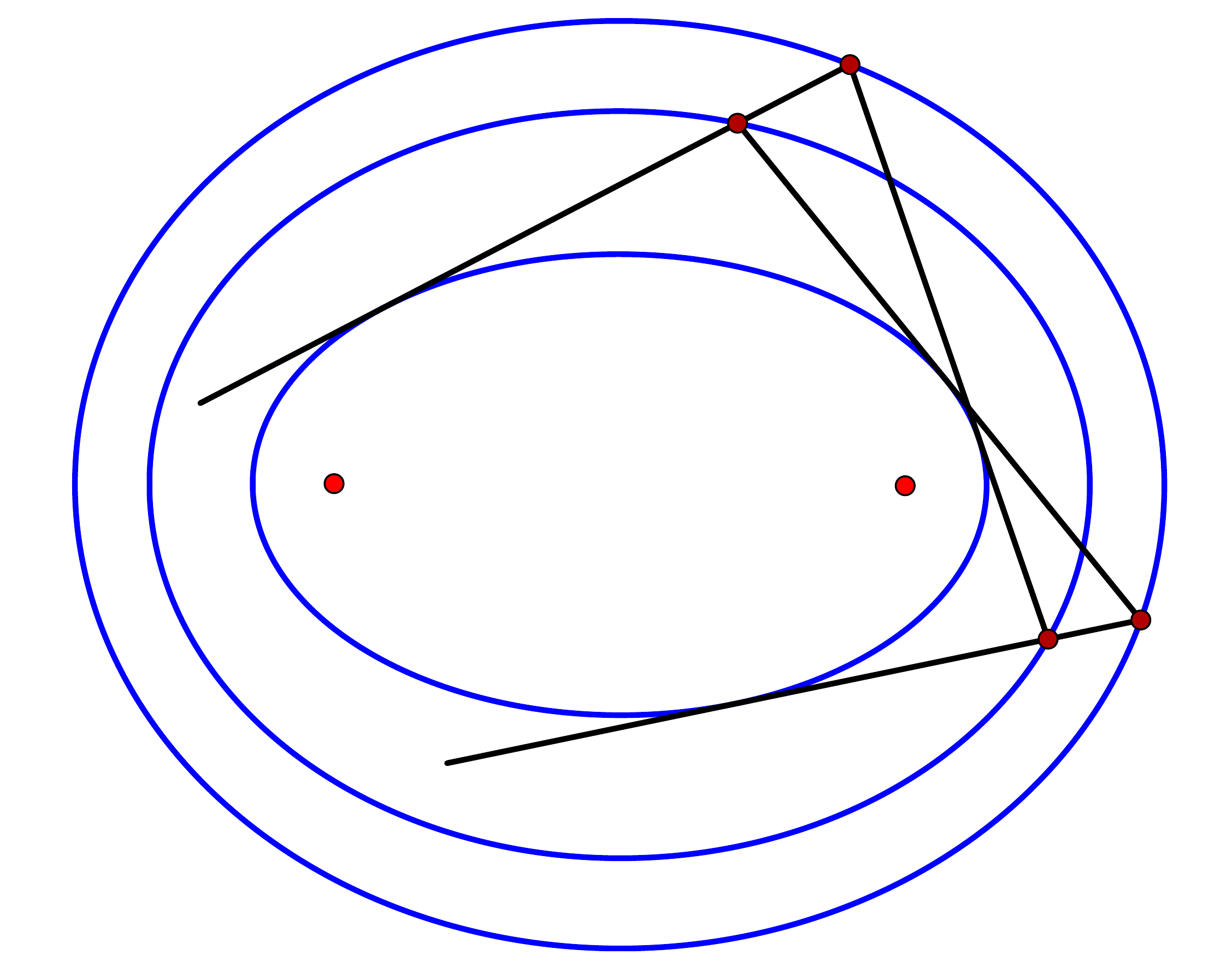}
\caption{The reflections from two confocal ellipses commute.}
\label{commute}
\end{figure}

Now we turn to the string diffeomorphisms on ellipses. Each ellipse $\g$ can be parameterized, $\g(t)$, such that, for every length of the string, the respective string diffeomorphism is a shift $t \mapsto t+c$, where the constant $c$ depends on the length of the string. In particular, these diffeomorphisms commute \cite[Corollary 4.6]{Tab05}, see Figure \ref{commute}.

We are led to the following definition.

\begin{definition} \label{def:Por}
{\rm
A closed smooth strictly convex curve $\g$ of length $\ell$ has the \emph{Poritsky property} if it has a parameterization $\g(t)$ such that for some $L > \ell$ and all lengths of the string in the interval $[\ell,L]$ the respective string diffeomorphisms are shifts $t \mapsto t+c$.
}
\end{definition}

The reason we call this the Poritsky property is that H. Poritsky proved that, in the Euclidean plane, a curve possessing this property is an ellipse \cite{Por} (his theorem is local: it concerns germs of plane curves, and concludes that the curves are conics). One of the results of \cite{Gl19} is an extension of the Poritsky theorem to the elliptic and hyperbolic planes and to outer billiards.

\begin{figure}[hbtp]
\centering
\includegraphics[height=1.7in]{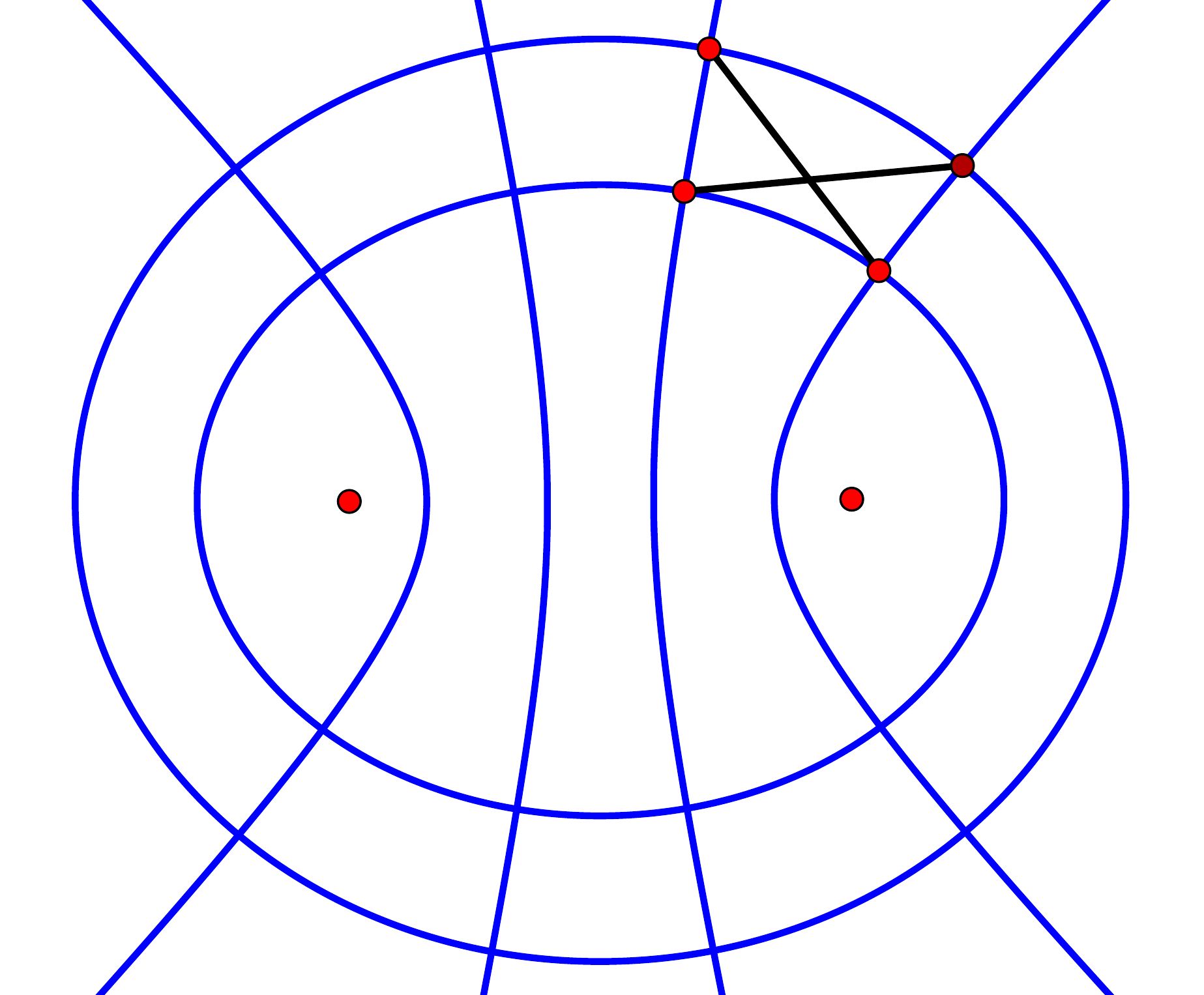}
\caption{Ivory's lemma.}
\label{Ivory}
\end{figure}

In addition to the ellipses with given foci, consider also the hyperbolas with the same foci. One obtains a net (or a $2$-web): through every point not on the axes there pass a unique ellipse and a unique hyperbola from the confocal family. A \emph{net quadrilateral} is a curvilinear quadrilateral whose sides lie on the net curves. The \emph{Ivory's lemma} says that {\it the diagonals of the curvilinear rectangles made by pairs of confocal ellipses and hyperbolas have equal length}, see Figure \ref{Ivory} and \cite[Theorem 1]{Iz17}.

This motivates the following definition.

\begin{definition}
A net of curves in the plane has the \emph{Ivory property} if the lengths of the diagonals in each net quadrilateral are equal.
\end{definition}

Theorems \ref{thm:main} and \ref{thm:LiouvEuc} stated on the next pages imply that if a closed curve in the plane possesses either Graves or Poritsky property, then it is an ellipse (in particular, this reproves Poritsky theorem). Similarly, the only nets possessing the Ivory property are the confocal nets and their degenerations.

The close link between the Graves, Poritsky, and Ivory properties becomes more apparent if one generalizes them to Riemannian surfaces.

\subsection{Riemannian surfaces and Liouville nets}
The string construction makes sense for any closed curve $\gamma$ without inflection points on a Riemannian surface. If points $A,B \in \g$ are sufficiently
close then the geodesics tangent to $\g$ at these points intersect at a point $C$ close to $\g$. We assume that the length of the string is sufficiently close to the length $\ell$ of the curve $\g$ so that the choice of the point $C$ is unique.
Under these assumptions, Definitions \ref{def:Graves} and \ref{def:Por} make sense in this more general setting.

Ellipses on the sphere or in the hyperbolic plane possess both the Graves and the Poritsky property \cite{Gl19}.

\begin{figure}[hbtp]
\centering
\includegraphics[height=1.7in]{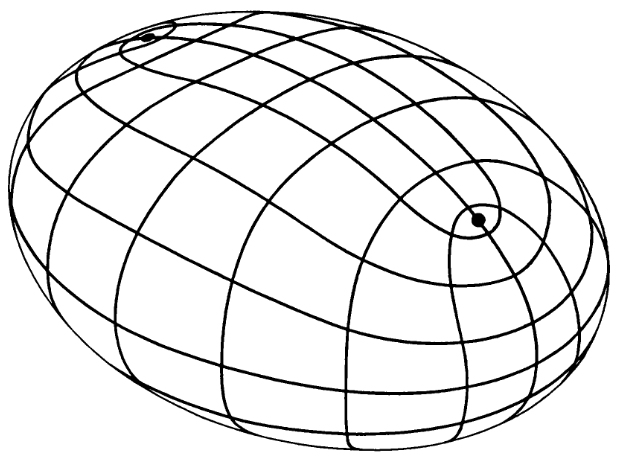}
\caption{Lines of curvature on an ellipsoid.}
\label{ellipsoid}
\end{figure}

A further example is provided by the metric on a 3-axial ellipsoid in $\R^3$. The lines of curvature are analogs of confocal conics, and the role of foci is played by the four umbilic points. See Figure \ref{ellipsoid}, borrowed from the classical book \cite{HC52}, for illustration. In the billiard table bounded by a line of curvature, the trajectory of a billiard ball is tangent to another line of curvature. Equivalently, the string construction around a line of curvature produces another line of curvature. Note that if $\g$ is a line of curvature, then there is a choice of the length of the string such that
the string construction yields the intersection of the ellipsoid with a coordinate plane. See \cite{ChSh89,Ves90,Ves91} concerning billiards on quadratic surfaces.

Note that the lines of curvature on an ellipsoid are the intersection curves with quadrics from the confocal family which contains this ellipsoid. The ellipses on the sphere and in the hyperbolic plane (in its hyperboloid model) are intersections with a certain family of quadratic cones obtained as a degeneration of a confocal family of quadrics.

Theorem \ref{thm:main} below states that the most general class of curves satisfying the Graves, Poritsky, and Ivory properties is given by \emph{Liouville nets}. Let us give the corresponding definitions.

\begin{definition} \label{def:Liouville}
A Riemannian metric in a 2-dimensional domain is called a \emph{Liouville metric} if there exist coordinates $(u,v)$ in which it is given  by the formula
\begin{equation} \label{eqn:Liouville}
(U_1(u) - V_1(v)) (U_2(u) du^2 + V_2(v) dv^2),
\end{equation}
where $U_i, V_i,\ i=1,2$ are smooth functions of one variable.
\end{definition}

An independent change of coordinates $\bar u = \bar u(u)$, $\bar v = \bar v(v)$ preserves the structure of the formula \eqref{eqn:Liouville}. This leads to the notion of a \emph{Liouville net}, which is a pair of orthogonal foliations such that for any choice of parametrization of leaves in each of them the resulting coordinates on the surface bring the Riemannian metric into the Liouville form.
An invariant (but quite complicated) local characterization of Riemannian surfaces possessing Liouville nets was given in \cite{Kru08}.

The Euclidean metric  is obviously a Liouville metric with Liouville net formed by lines 
parallel  to orthogonal coordinate axes. It also admits  other Liouville nets  
formed by confocal conics. Indeed,
a family of confocal conics is given by the formula
$$
\frac{x^2}{a+\lambda} + \frac{y^2}{b+\lambda} =1,\ a>b>0,
$$
where $\lambda$ is the parameter. The parameter values corresponding to two conics passing through a given point are the \emph{elliptic coordinates} of the point. In the elliptic coordinates, the Euclidean metric $dx^2+dy^2$ has the following form:
\begin{equation} \label{ellcoord}
(\lambda - \mu) \left( \frac{d\lambda^2}{4(a+\lambda)(b+\lambda)} - \frac{d\mu^2}{4(a+\mu)(b+\mu)} \right).
\end{equation}

The ellipsoidal coordinates $(\lambda, \mu, \nu)$ in $\R^3$ restricted to the ellipsoid $\nu = \const$ describe a Liouville metric whose coordinate curves are the lines of curvature on the ellipsoid. Liouville introduced the metrics \eqref{eqn:Liouville} as an immediate generalization of the ellipsoidal coordinates: geodesic equation can be solved for Liouville metrics by the Jacobi separation method. Higher-dimensional analogs of Liouville metrics/nets are St\"ackel metrics/nets.

See \cite{Kiy97} for general material about Liouville metrics and \cite{Iz17,PoTo03,PoTo08,PoTo11} for Liouville billiards.

Now we formulate our main result.

Consider a closed annulus ${\mathcal A}$ equipped with a Riemannian metric and foliated by smooth closed geodesically convex curves. (By smooth we always mean $C^\infty$-smooth.) Call this foliation ${\mathcal F}_1$, and let ${\mathcal F}_2$ be the foliation formed by the curves orthogonal to the leaves of ${\mathcal F}_1$.

\begin{theorem} \label{thm:main}
The following four properties are equivalent:\\
(i) The foliation ${\mathcal F}_1$ has the Graves property;\\
(ii) The inner boundary curve of ${\mathcal A}$ has the Poritsky property;\\
(iii) The foliations ${\mathcal F}_1$ and ${\mathcal F}_2$ form a Liouville net;\\
(iv) The net $({\mathcal F}_1,{\mathcal F}_2)$ in ${\mathcal A}$ has the Ivory property.
\end{theorem}

The conjecture, attributed to Birkhoff, states that {\it if the interior neighborhood of a smooth strictly convex boundary curve of a billiard table in the Euclidean plane is foliated by caustics, then this curve is an ellipse (or a circle).} A substantial progress has been made recently toward the proof of this conjecture; see the surveys \cite{BM18,KS18} and the references therein.

We finish this introduction with a version of Birkhoff conjecture.

\begin{conjecture} \label{conj1}
Consider an annulus equipped with a Riemannian metric in which one of the components $\Gamma$ of the boundary is strictly convex, and consider the billiard system near this component. If a neighborhood of $\Gamma$ is foliated by caustics, then the metric near $\Gamma$ is Liouville, and $\Gamma$ is one of the coordinate lines.
\end{conjecture}

Conjecture \ref{conj1} implies the classical Birkhoff Conjecture due to the following classical result.

\begin{theorem}[Wei24, Bla28]
\label{thm:LiouvEuc}
Any Liouville net in a domain $\Omega \subset \R^2$ is of one of the following types.
\begin{enumerate}
\item
Net of confocal ellipses and hyperbolas.
\item
Net of confocal and coaxial parabolas.
\item
Net of concentric circles and their radial lines.
\item
Orthogonal net of lines.
\end{enumerate}
\end{theorem}

Weihnacht's result is more general: he deals with Liouville nets in dimension $2$ and with St\"ackel nets in dimension $3$. Blaschke discusses only the $3$-dimensional case using the equivalence of the Liouville and Ivory properties. In the Appendix we give a proof of Theorem \ref{thm:LiouvEuc} based on Blaschke's ideas.

\subsection{Local version of the main result}

Let $\mcd$ be a topological disc equipped with a smooth 
Riemannian metric. Let  $\gamma\subset \mcd$ be a germ of curve with 
positive geodesic curvature at a point $O\in \mcd$. For every  two given points $A,B\in\gamma$,  we will denote by $C_{AB}$ the  point of intersection 
of the geodesics tangent to $\gamma$ at $A$ and 
$B$. We consider the situation when $A$ and $B$ are close to each other, and then we can, and will, choose the above $C_{AB}$ close to $A$ and $B$ in a unique way. Set 
$$\la(A,B):= \text{ the length of the arc } AB \text{ of the curve } \gamma, $$
\begin{equation}L(A,B):=|AC_{AB}|+|BC_{AB}|-\la(A,B).\label{lab}\end{equation}
Here, for $X,Y\in\mcd$ close enough to a given point $O$, 
 we denote by $|XY|$ the length of the geodesic segment connecting $X$ and $Y$.
 
\begin{definition} \label{dstring} Under the above conditions, 
for every $p\in\rr_+$ small enough, 
the subset  
$$\Gamma_p:=\{ C_{AB} \ | \  L(A,B)=p\}\subset\mcd$$
is  called the {\it $p$-th string construction}; we set $\Gamma_0=\gamma$
\end{definition}

\begin{remark} In the case when $\gamma$ is a closed curve in a Riemannian surface, 
for every $O\in\gamma$, the family of curves given by the global string construction 
(Subsection 1.1) coincides locally near $O$ with the above family $\Gamma_p$ 
given by the local string construction.
\end{remark}
\begin{remark} For every $p>0$ small enough,  $\Gamma_p$ is a well-defined smooth curve.  
The curve $\gamma$ is a caustic for the billiard in the curve $\Gamma_p$, as in the 
global string  construction. 
Let $U$ be a small neighborhood of the base point $O$, and let $\mca$ be the connected 
component of the complement $U\setminus\gamma$ lying on the 
concave side from the curve $\gamma$.  As in the closed curve case,  
the curves $\Gamma_p$ with small $p\geq0$ form a smooth foliation $\mathcal F_1$ 
on $\mca$  by smooth 
curves. \end{remark}

The definitions of Graves and Poritsky properties of the above foliation on $\mca$ 
by curves $\Gamma_p$ are the same as in the previous subsection.  

\begin{theorem} \label{thm:main2} Let $\gamma$ be a germ of curve in a Riemannian 
topological disc. Let $\mca$ be as above. Let $\mathcal F_1$ be the foliation 
of the domain $\mca$ by  curves $\Gamma_p$  given by the above local string construction. 
Let $\mathcal F_2$ be the foliation orthogonal to $\mathcal F_1$. 
The following four properties are equivalent:\\
(i) The foliation ${\mathcal F}_1$ has the Graves property;\\
(ii) The curve $\g$ has the Poritsky property;\\
(iii) The metric in ${\mathcal A}$ is Liouville, and the leaves of the foliations ${\mathcal F}_1$ and ${\mathcal F}_2$ are the coordinate lines;\\
(iv) The net $({\mathcal F}_1,{\mathcal F}_2)$ in ${\mathcal A}$ has the Ivory property.
\end{theorem}

\subsection{Plan of the paper and the proof} 

The proof of Theorem \ref{thm:main} given below is done by  local arguments and 
 remains valid in the local case: for Theorem  \ref{thm:main2}. 
 
Implication (i)$\Rightarrow$(ii) is proved in Section 2. 

Implication (ii)$\Rightarrow$(iii) is proved in Section 3. 

Equivalence of statements (iii) and (iv) is proved in Section 4. 

Equivalence (iii)$\Leftrightarrow$(i) is explained in Section 5. 

\section{Graves implies Poritsky} \label{sect:GtoP}

In this section we shall show that the Graves property implies the Poritsky property, that is, we establish the implication (i) $\Rightarrow$ (ii) in Theorem \ref{thm:main}. The argument is not new and, in various forms, it is described in the literature; see, e.g., \cite{Tab95}.
We break the argument into a number of lemmas, and we only sketch the proofs, referring for details to the literature.

\begin{lemma} \label{lem:spgeo} Let $\mcd$ be a Riemannian surface. For every 
point $O\in\mcd$ and every strictly convex neighborhood $U=U(O)\subset\mcd$ 
with smooth boundary, the space of oriented geodesics in $U$ is a two-dimensional 
manifold. It is topologically a cylinder  $S^1 \times (0,1)$.
\end{lemma}

\proof
The convexity assumption   on $U$ implies that every oriented geodesic in $U$ intersects its boundary at two points: the departure point and the arrival point. 
To each oriented geodesic in $U$ we associate the pair $(x,v)$: $x\in\partial U$ being the 
departure point and $v$  its orienting unit tangent  vector at $x$. For a given $x\in\partial U$ the space of all the departure vectors $v\in T_x\mcd$ 
 is the semicircle  of unit vectors 
directed inside $U$.  The above correspondence 
is bijective. Its image is a cylinder: 
the product of the boundary $\partial U$ and  semicircle. 
\proofend

Denote the above space of oriented geodesics in $U$ by $\mathcal G={\mathcal G}_U$. 

\begin{lemma} \label{lem:symp1} (Melrose construction, see \cite{Me76}). 
The space ${\mathcal G}$ has an area form $\omega$ obtained by the symplectic reduction from the canonical symplectic form in the cotangent bundle $T^* {\mathcal D}$. 
The area forms coming from different intersecting neighborhoods coincide. 
\end{lemma}

\proof
 The construction described below works in any dimension.
  
The cotangent bundle $T^* {\mathcal D}$ carries the canonical symplectic structure 
$\Omega=dp \wedge dq$, where $q$ are the positions and $p$ are the momenta. The Riemannian metric makes it possible to talk about the norms of covectors.
The geodesic flow is the Hamiltonian flow with the energy function $|p|^2/2$. Identifying the tangent and cotangent bundles by the Riemannian metric, we obtain the geodesic flow on the tangent bundle $T{\mathcal D}$.

Consider the unit energy hypersurface $S \subset T^* {\mathcal D}$. The restriction of the symplectic form $\Omega$ on $S$ has 1-dimensional kernel, and $S$ is foliated by the integral curves of this line field, the characteristic curves. Each characteristic curve is the trajectory of the Hamiltonian vector field of the energy function, and it can be viewed as an arc length  parameterized geodesic. 

Let $\pi:T\mcd\to\mcd$ denote the standard projection. 
The quotient space of $S\cap\pi^{-1}(U)$ by the characteristic foliation is the space of oriented non-parameterized geodesics in $U$, and since one factorizes by the kernel of the restriction of the symplectic structure $\Omega$ on $S$, we obtain an induced symplectic structure on  ${\mathcal G}$. The last statement of the lemma follows immediately 
from definition. 
\proofend

Let $\delta$ be a geodesically convex simple closed curve in ${\mathcal D}$ 
(or a germ of curve with positive geodesic curvature at a point $O\in\mcd$). 
Consider the billiard inside $\delta$. 
The respective billiard ball map $F$ acts on the space of 
oriented geodesics that intersect $\delta$. 
(In the case of germ it acts on the oriented geodesics intersecting $\delta$ 
that are close to the geodesic 
tangent to $\delta$ at $O$.) Abusing notation, we denote the restriction of the symplectic form $\omega$ on the latter space of geodesics by the same letter. 

\begin{lemma} \label{lem:symp2}
The billiard ball map $F$ preserves the form $\omega$.
\end{lemma}

\proof
Once again, the argument works in any dimension, and it is an adaptation of a more general result of Melrose \cite {Me76}.

Assign to an oriented geodesic intersecting $\delta$ its first intersection point $A$ with $\delta$ and its unit tangent vector $v$ at this point. This gives an identification of this space of geodesics with the space of unit tangent vectors with foot point on $\delta$ having the interior direction. Projecting the tangent vector $v$ to the tangent space $T_A \delta$ identifies the space of geodesics with the unit disk subbundle of $T \delta$, identified with $T^* \delta$ by the metric. The latter has its own canonical symplectic structure ``$dp \wedge dq$".  One can prove that this structure coincides with $\omega$; see \cite{Tab95}.

The billiard ball map is a composition of two transformations. The first one 
takes the a unit vector $v\in T_A\mcd$ to a unit vector $w\in T_B\mcd$, where $B$ is the second intersection point of the oriented geodesic tangent to $v$ with $\delta$ and $w$ 
is its tangent vector at $B$. The second transformation 
reflects the vector $w$ from the tangent space $T_B \delta$. The first transformation is seen to preserve the form $\omega$, constructed in Lemma \ref{lem:symp1}. The reflection of a unit tangent vector in tangent space $T_B \delta$ does not change its projection to  $T_B \delta$, hence it also preserves $\omega$. This proves the lemma.
\proofend

Let us reiterate: for every billiard table, the billiard ball map  preserves the same area form on the space of oriented geodesics; this form  depends only on the Riemannian metric.

\begin{remark}
{\rm Another way to construct an invariant symplectic form of the billiard ball map is by considering the billiard inside $\delta$ as a discrete Lagrangian system. The geodesic $AB$ reflects to the geodesic $BC$, where $A,B,C \in \delta$, if and only if the sum of the Riemannian distances 
$|AB|+|BC|$ has a critical point at $B$ as a function of $B \in \delta$. 
That is, 
$$
L(A,B)=|AB|: \delta \times \delta \to \R
$$
is the Lagrangian function. 
Then the 2-form 
$$
\frac{\partial^2 L(A,B)}{\partial A\ \partial B}\ dA \wedge dB 
$$
is invariant under the billiard ball map $(A,B) \mapsto (B,C)$, see \cite[chapter 1, section 1]{Ves911}. Here the derivatives and the differentials are taken with respect to the length 
parameters of the points $A$ and $B$. 

However, in this approach, it is not clear that the invariant 2-form is the same for different billiard curves $\delta$. 
}
\end{remark}

Let $\g_1$ be a leaf of the foliation ${\mathcal F}_1$. Consider the billiard system inside $\g_1$. The leaves of ${\mathcal F}_1$ that lie between $\g$ and $\g_1$ are caustics of this billiard. Consider the part of the phase space that is the union of the respective invariant curves under the billiard ball map: it looks like the lower part of Figure \ref{phase}, below the singular, eye-shaped, leaf. The billiard ball map gives a transformation of each of these invariant curves. 

\begin{lemma} \label{lem:para}
Each invariant curve can be parameterized so that the billiard ball transformation is a shift $t \mapsto t+c$.
\end{lemma} 

\proof
Choose a function $H$  whose level curves are the invariant curves, and consider its Hamiltonian vector field sgrad $H$ with respect to the area form $\omega$. This vector field is tangent to the invariant curves, and the desired coordinate $t$ is the one in which sgrad $H$ is a constant vector field $d/dt$. Changing $H$ scales the coordinate $t$ on each invariant curve. One can normalize the domain of the parameter $t$ to be the unit circle, and this  fixes $t$  up to an additive constant. In other words, the 1-form $dt$ is well defined on each invariant curve. 

The billiard ball map preserves $\omega$ and the invariant curves, hence it preserves $dt$, that is, the map is a shift $t \mapsto t+c$.
\proofend

As an aside, we note two geometric consequences. The first is a version of the Poncelet porism: if a billiard trajectory in $\g_1$ closes up after a number of reflections, then all the trajectories with the same caustic close up after the same number of reflections. 

The second consequence concerns the billiard systems inside two leaves of the foliation ${\mathcal F}_1$, curves $\g_1$ and $\g_2$. These billiards share the caustics, and the  invariant measure $dt$ on the respective invariant curves depends only on the metric and the foliation by invariant curves, that is, this invariant measure is the same for both billiards. Since shifts commute, it follows that the two billiard ball transformations commute as well; 
this argument is presented in \cite[pp. 58, 59, corollary 4.6]{Tab05}. 

Thus the situation is exactly the same as for ellipses in the Euclidean plane.

Finally we show that the curve $\g$ has the Poritsky property. Indeed, $\g$ is a caustic for the billiard system inside each leaf of the foliation ${\mathcal F}_1$, hence it carries a parameter $t$ constructed in Lemma \ref{lem:para}. As we already mentioned, this parameter is the same for every choice of a leaf. Invoking the string construction, which  recovers the billiard table from its caustic, we see that $\g$ possesses the Poritsky property. 

\section{Poritsky implies Liouville} \label{sect:PtoL}

Let $\g$ possess the Poritsky property, and let $\g(t)$ be the respective parameterization. 

\begin{figure}[hbtp]
\centering
\includegraphics[height=1.7in]{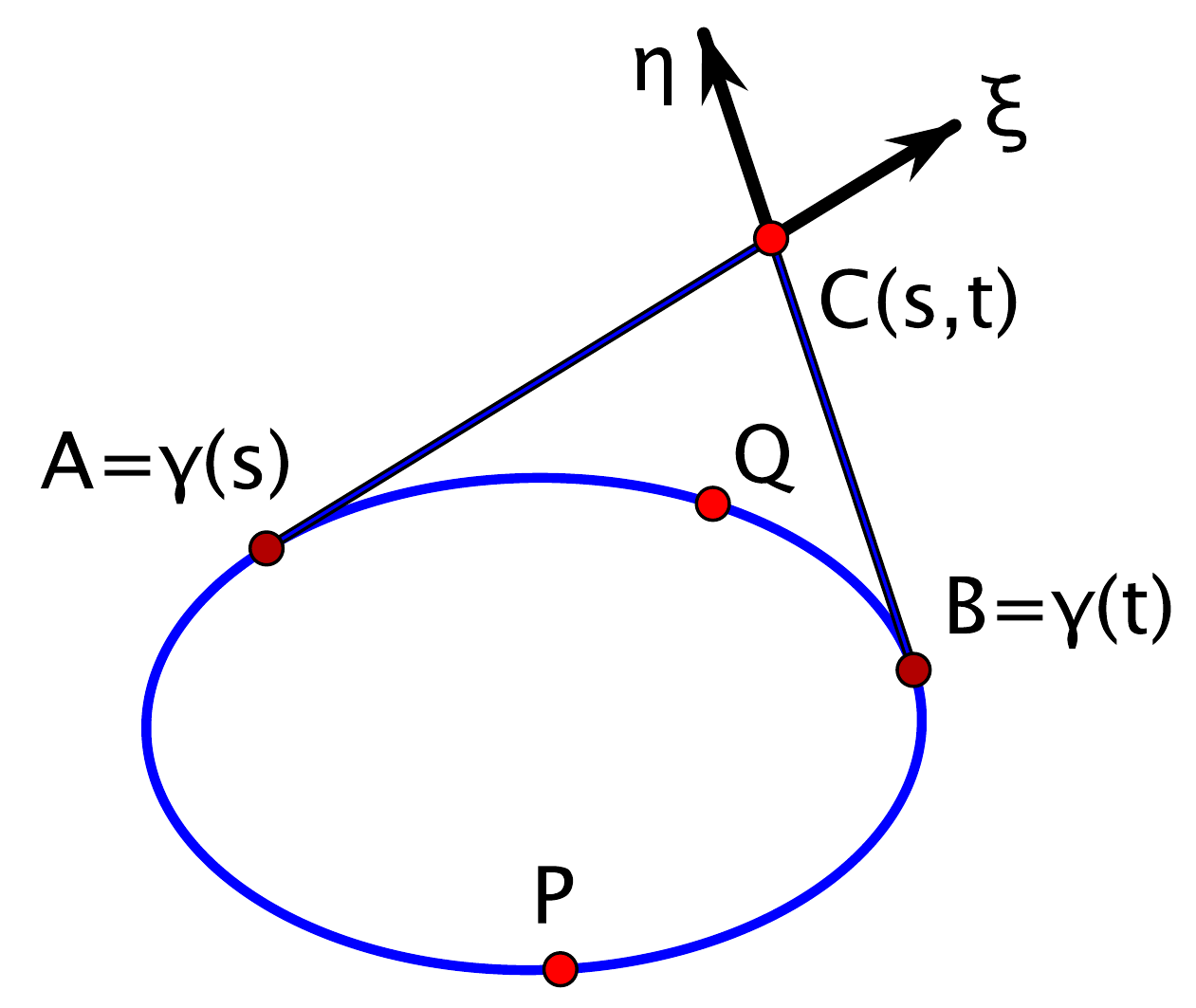}
\caption{Coordinates in the exterior of the curve $\g$.}
\label{coord}
\end{figure}

Consider Figure \ref{coord}: we can use $(s,t)$ as coordinates of point $C$. We also consider another coordinate system
$$
x=\frac{s+t}{2},\ y= \frac{t-s}{2}
$$
(these formulas make sense locally; globally, they are well-defined up to adding a constant to the $x$-coordinate).

\begin{lemma} \label{lem:extcoord} The metric in the domain ${\mathcal A}$ admits orthogonal coordinates $(x,y)$ in which the diagonals $x\pm y=const$ are geodesics. 
\end{lemma}

\proof
Let $\xi$ and $\eta$ be the unit tangent vectors along the geodesics $AC$ and $BC$ in Figure \ref{coord}. These geodesics are given by the equations $s=$ const and $t=$ const. Hence $ds(\xi)=dt(\eta)=0$.

Each curve given by the string construction has the equation $t-s=$ const, and the billiard property implies that the vector $\xi-\eta$ is tangent to this curve. Hence $(dt-ds)(\xi-\eta)=0$. This implies that $dt(\xi) + ds(\eta)=0$, and therefore $(dt+ds)(\xi+\eta)=0$. 

Since the vectors $\xi$ and $\eta$ are unit,  $\xi+\eta$ is orthogonal to $\xi - \eta$. But the vector $\xi+\eta$ is tangent to the level curve of $s+t$. Therefore this level curve is orthogonal to the level curve of $t-s$. 

Changing to the $(x,y)$-coordinates, we see that the coordinate lines $x=$ const and $y=$ const are orthogonal. And the geodesics $AC$ and $BC$ are the diagonals $x-y=$ const and $x+y=$ const, as claimed. 
\proofend

\begin{remark} 
{\rm In Figure \ref{coord}, consider two functions, $\varphi(C)= |AC| + \stackrel{\smile}{|AP|}$ and $\psi(C)= |BC| + \stackrel{\smile}{|BP|}$, the distances from point $C$ to point $P$ 
going around the obstacle $\g$; here $\stackrel{\smile}{|AP|}$ and $\stackrel{\smile}{|BP|}$ 
are the  lengths of the corresponding lower arcs of the curve $\gamma$. 
The vectors $\xi$ and $\eta$ are the gradients of these functions. Each curve given 
by the string construction is a level curve of the function $\varphi+\psi$. 
This implies  that the orthogonal curves are the level curves of the function $\varphi-\psi$. See \cite{LT07,Iz17} for an argument in the Euclidean case which applies in the more general Riemannian case as well.

We add that that one does not need to consider the ``invisible part" of the curve $\g$: one can redefine the functions $\varphi$ and $\psi$ by taking the difference $|AC| - \stackrel{\smile}{|AQ|}$ and $|BC| - \stackrel{\smile}{|BQ|}$. This allows to adapt our arguments for the proof of Theorem \ref{thm:main2}. 
}
\end{remark}

We will now prove the main result of this section.

\begin{theorem} \label{thm:2} 
Let a Riemannian metric, written in orthogonal coordinates $(x,y)$, have the property that 
the diagonals $x\pm y=const$ are geodesics. Then this metric is Liouville. 
\end{theorem}

\proof
Let the metric be 
\begin{equation} 
a(x,y)dx^2+b(x,y)dy^2.\label{metric}
\end{equation}

The constant speed parametrized geodesics  are the extremals of the functional 
$$
\int_0^1L(x,y,\dot x,\dot y)dt, \ L(x,y,\dot x,\dot y)=a(x,y)\dot x^2+b(x,y)\dot y^2.
$$
The corresponding system of Euler--Lagrange equations 
$$
\frac{d}{dt}\frac{\partial L}{\partial \dot x}=
\frac{\partial L}{\partial x}, \ \frac{d}{dt}\frac{\partial L}{\partial \dot y}=
\frac{\partial L}{\partial y}
$$ 
takes the following form:
\begin{equation}\begin{cases} 
2a\ddot x+2\dot x^2\frac{\partial a}{\partial x}+2\dot x\dot y
\frac{\partial a}{\partial y}= \dot x^2\frac{\partial a}{\partial x}+\dot y^2\frac{db}{dx}\\
2b\ddot y+2\dot y^2\frac{\partial b}{\partial y}+2\dot x\dot y
\frac{\partial b}{\partial x}= \dot y^2\frac{\partial b}{\partial y}+\dot x^2\frac{da}{dy}.
\end{cases}\label{el1}
\end{equation}

Let us write down equations (\ref{el1}) on the diagonals $\pm x-y=$ const, oriented by the coordinate $x$, 
which are known to be geodesics.  Along these diagonals one has 
$$
\dot x=\pm\dot y, \ \ddot x=\pm\ddot y.
$$
Substituting the latter formulas to (\ref{el1}), we get: 
\begin{equation}\begin{cases} 
2a\ddot x=\dot x^2(\frac{\db}{\dx}-\frac{\da}{\dx}\mp2\frac{\da}{\dy})\\
\pm 2b\ddot x=\dot x^2(\frac{\da}{\dy}-\frac{\db}{\dy}\mp 2\frac{\db}{\dx}).
\end{cases}\label{el2}
\end{equation}

Let us multiply the equations in (\ref{el2}) by $b$ and $a$, respectively. 
In the case of "$+$" (the diagonals $y-x=$ const), taking the difference of thus obtained equations,   yields 
\begin{equation} 
b\left(\frac{\db}{\dx}-\frac{\da}{\dx}-2\frac{\da}{\dy}\right)=
a\left(\frac{\da}{\dy}-\frac{\db}{\dy}-2\frac{\db}{\dx}\right).
\label{+}\end{equation}
In the case of "$-$" (the diagonals $x+y=$ const), taking their sum, yields 
\begin{equation} 
b\left(\frac{\db}{\dx}-\frac{\da}{\dx}+2\frac{\da}{\dy}\right)=
-a\left(\frac{\da}{\dy}-\frac{\db}{\dy}+2\frac{\db}{\dx}\right).
\label{-}\end{equation}
Taking the sum and the difference of equations (\ref{+}) and (\ref{-}) yields the system 
\begin{equation}\begin{cases}
b(\frac{\db}{\dx}-\frac{\da}{\dx})=-2a\frac{\db}{\dx}\\
2b\frac{\da}{\dy}=a(\frac{\db}{\dy}-\frac{\da}{\dy}).
\end{cases}\label{pm}\end{equation}
System (\ref{pm}) can be rewritten as 
\begin{equation} \begin{cases}
\frac{\da}{\dx}=(1+2\frac{a}b)\frac{\db}{\dx}\\
\frac{\db}{\dy}=(1+2\frac ba)\frac{\da}{\dy}.
\end{cases}
\label{pm2}\end{equation}

The first equation in (\ref{pm2}) says that, along the lines parallel to the $x$-axis, the implicit function $a(b)$ 
satisfies the differential equation
\begin{equation}\frac{da}{db}=1+2\frac ab.\label{diff1}\end{equation}
Each solution of equation (\ref{diff1}) is of the form $a(b)=c_1b^2-b$. Therefore,
\begin{equation} 
a=c_1(y)b^2-b \label{axb}.
\end{equation}
Similarly the second equation in (\ref{pm2}) says that, along the lines parallel to the $y$-axis, the implicit 
function $b(a)$ satisfies a similar differential equation, and we get that 
\begin{equation} 
b=c_2(x)a^2-a\label{bxa}.
\end{equation}

Combining formulas (\ref{axb}) and (\ref{bxa}) we obtain 
$$
b=c_2(x)a^2-a=c_2(x)(c_1(y)b^2-b)^2-c_1(y)b^2+b.
$$
The latter is equivalent to the equality 
$$
c_2(x)(c_1(y)b-1)^2=c_1(y).
$$
Thus we obtain that 
$$
a(x,y)= \frac{f(x)+g(y)}{f^2(x)g(y)}, \ b(x,y)=\frac{f(x)+g(y)}{f(x)g^2(y)},
$$
where
$$
f(x)=(c_2(x))^{\frac12}, g(y)=(c_1(y))^{\frac12}.
$$

Finally, make a coordinate change $u(x),v(y)$ given by the formula
$$
\frac{dx}{f(x)^{\frac12}} = du,\ \frac{dy}{g(y)^{\frac12}} = dv.
$$
Then the metric (\ref{metric}) becomes
$$
\left( \frac{1}{f(x)} + \frac{1}{g(y)} \right) (du^2+dv^2) = (U(u)-V(v))(du^2+dv^2),
$$
which has the Liouville form (\ref{eqn:Liouville}). 
This finished the proof.
\proofend

Lemma \ref{lem:extcoord} and Theorem \ref{thm:2} combined show that the Poritsky property of the curve $\g$ implies that the metric in the annulus ${\mathcal A}$ is Liouville.

\section{Ivory is equivalent to Liouville} \label{sect:IandL} 

The theorem that the Liouville nets are characterized by the Ivory property is due to Blaschke and Zwirner \cite{Bla28,Zw27}; see also \cite{Thimm78}. In \cite{Iz17}, the 
second and the third authors of the present paper provided a streamlined account of the implication Liouville $\Rightarrow$ Ivory. Here we give a proof of the converse implication based on Blaschke's arguments in \cite[\S 56]{BlaschkeBook}. See also \cite{Bla28, Thimm78} for the higher-dimensional analog: Ivory property holds only in St\"ackel nets.

\begin{theorem}
\label{thm:IvoryLiouville}
If a Riemannian metric $g$ satisfies the Ivory property, then it is Liouville:
\[
g = (U_1 - V_1)(U_2 du^2 + V_2 dv^2),
\]
where $U_1, U_2$ are functions of $u$, and $V_1, V_2$ are functions of $v$.
\end{theorem}

We start with a number of lemmas.

\begin{lemma}
\label{lem:IvoryOrth}
In a metric with Ivory property, the coordinate lines are orthogonal.
\end{lemma}

\proof
Let $(u,v)$ be a point sufficiently close to $(u_0,v_0)$, and let 
$a\,du^2 + 2b\,du\,dv + c\,du^2$ be the metric tensor at $(u,v)$. Let $\varepsilon > 0$ be 
sufficiently small, so that the coordinate rectangle $[u,u+\varepsilon] \times [v,v+\varepsilon]$ 
has diagonals of equal length. 
As $\varepsilon$ tends to zero, the diagonal lengths satisfy
\[
L_+ \sim \varepsilon \sqrt{a+2b+c}, \quad L_- \sim \varepsilon \sqrt{a-2b+c},
\]
where $L_+$ is the length of the diagonal containing $(u,v)$, and $L_-$ is the length of the diagonal containing $(u,v+\varepsilon)$. The equality $L_+ = L_-$ implies $b = 0$, and we are done.
\proofend

\begin{lemma}
\label{lem:IvoryLem1}
Let $[u_1,u_2] \times [v_1,v_2]$ be a coordinate rectangle satisfying the Ivory property,
and let $\gamma_+ \colon [0,L] \to \R^2$, $\gamma_- \colon [0,L] \to \R^2$
be its (arc length parametrized) geodesic diagonals such that
\begin{align*}
\gamma_+(0) = (u_1,v_1), \quad &\gamma_+(L) = (u_2,v_2),\\
\gamma_-(0) = (u_1,v_2), \quad &\gamma_-(L) = (u_2,v_1).
\end{align*}
If all sufficiently close rectangles also satisfy the Ivory property, then one has
\[
\left\langle \dot{\gamma_+}(L), \frac{\partial}{\partial u} \right\rangle = \left\langle \dot{\gamma_-}(L), \frac{\partial}{\partial u} \right\rangle, \quad
\left\langle \dot{\gamma_+}(L), \frac{\partial}{\partial v} \right\rangle = -\left\langle \dot{\gamma_-}(0), \frac{\partial}{\partial v} \right\rangle.
\]
\end{lemma}

\proof
Deform the rectangle by moving its $u=u_2$ side.
Let $\gamma_+^\varepsilon$ be the geodesic from $(u_1,v_1)$ to $(u_2+\varepsilon,v_2)$,
and $\gamma_-^\varepsilon$ be the geodesic from $(u_1,v_2)$ to $(u_2+\varepsilon,v_1)$. 
Let $L(\gamma^{\varepsilon}_{\pm})$ denote the length of the 
curve $\gamma^{\varepsilon}_{\pm}$. 

\begin{figure}[hbtp]
\centering

\includegraphics[height=1.7in]{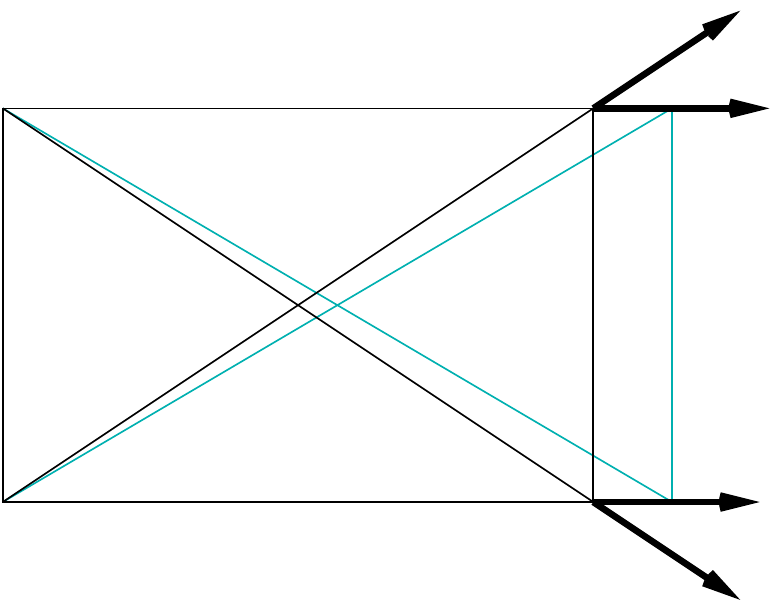}
\caption{To the proof of Lemma \ref{lem:IvoryLem1}.}
\label{fig:Lem1}
\end{figure}

Then, by the first variation of arc length formula, one has
\[
\left.\frac{d}{d\varepsilon}\right|_{\varepsilon = 0} L(\gamma_+^\varepsilon) = \left\langle \dot{\gamma_+}(L), \frac{\partial}{\partial u} \right\rangle, \quad
\left.\frac{d}{d\varepsilon}\right|_{\varepsilon = 0} L(\gamma_-^\varepsilon) = \left\langle \dot{\gamma_-}(L), \frac{\partial}{\partial u} \right\rangle.
\]
By the Ivory property, one has $L(\gamma_+^\varepsilon) = L(\gamma_-^\varepsilon)$ for all sufficiently small $\varepsilon$.
This implies the first equation of the Lemma. See Figure \ref{fig:Lem1}.
The second equation is proved similarly by moving the $v = v_2$ side of the rectangle. 
\proofend

Lemma \ref{lem:IvoryLem1} can be reformulated as follows:
at each of the corners of the coordinate rectangle, take the velocity vector of the geodesic diagonal ending at this corner
and turn it into a covector via the Riemannian metric.
Then the $du$-components of the covectors at the corners with the same $u$-coordinate coincide, up to the sign,
and the $dv$-components of the covectors at the corners with the same $v$-coordinate coincide, up to the sign.

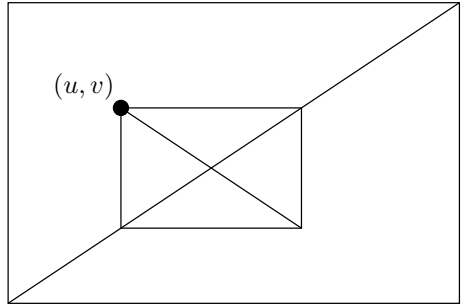
\begin{figure}[hbtp]
\centering
\begin{picture}(0,0)%
\includegraphics{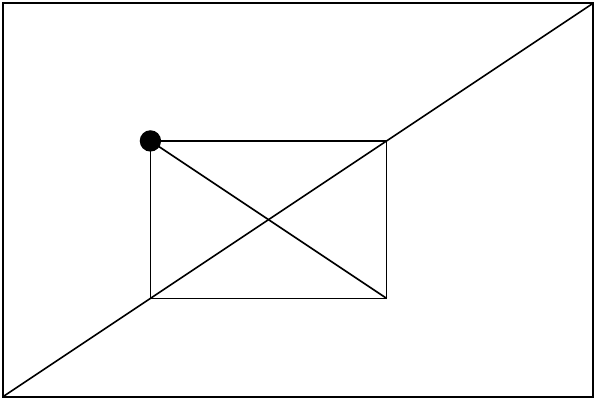}%
\end{picture}%
\setlength{\unitlength}{4144sp}%
\begingroup\makeatletter\ifx\SetFigFont\undefined%
\gdef\SetFigFont#1#2#3#4#5{%
  \reset@font\fontsize{#1}{#2pt}%
  \fontfamily{#3}\fontseries{#4}\fontshape{#5}%
  \selectfont}%
\fi\endgroup%
\begin{picture}(2724,1824)(-11,-973)
\put(271,299){\makebox(0,0)[lb]{\smash{{\SetFigFont{10}{12.0}{\rmdefault}{\mddefault}{\updefault}{$(u,v)$}%
}}}}
\end{picture}%
\caption{To the definition and properties of $\eta_\pm$.}
\label{fig:Lem2}
\end{figure}

Let $[u_1,u_2] \times [v_1,v_2]$ be a coordinate rectangle such that all coordinate rectangles inside of it satisfy the Ivory property. Let $\gamma_+$ be its geodesic 
diagonal, as in Lemma \ref{lem:IvoryLem1}. 
Define two differential $1$-forms in the rectangle as follows.
Through any point $(u,v)$ draw coordinate lines until they intersect the geodesic diagonal $\gamma_+$. Let $(u,v')$ and $(u',v)$ denote the corresponding intersection 
points. 
Take the velocity vector of $\gamma_+$ at the point $(u,v')$, and denote by $\phi(u)$ the $du$-component of the corresponding covector.
Similarly, denote by $\psi(v)$ the $dv$-component of the covector dual to the velocity vector of $\gamma_+$ at $(u',v)$.
Then put
\[
\eta_+(u,v) = \phi(u) du + \psi(v) dv, \quad \eta_-(u,v) = \phi(u) du - \psi(v) dv.
\]
See Figure \ref{fig:Lem2}.

\begin{lemma}
\label{lem:IvoryLem2}
One has $\|\eta_\pm\| = 1$.
\end{lemma}

\proof
Consider the rectangle with corners $(u,v)$, $(u,v')$, $(u',v)$, and $(u',v')$. 
Draw a geodesic diagonal from $(u,v)$ to $(u',v')$.
By Lemma \ref{lem:IvoryLem1} and by definition of $\eta_-$, the velocity vector of this geodesic at $(u,v')$ is the image of $\eta_-$
under the isomorphism induced by the Riemannian metric.
This implies $\|\eta_-\| = 1$. Due to the orthogonality of the coordinate system, one also has $\|\eta_+\| = 1$.
\proofend

\begin{corollary}
Let $f(u)$ and $g(v)$ be functions such that $\phi(u) du = df$ and $\psi(v) dv = dg$.
The integral curves of the gradient fields $\nabla(f+g)$ and $\nabla(f-g)$ are geodesics.
\end{corollary}

\proof
The gradient of $f \pm g$ has norm $1$ everywhere.
By a well-known theorem, the integral curves of such a gradient field are geodesics.
\proofend

%

\noindent \textbf{Proof of Theorem \ref{thm:IvoryLiouville}.}
By Lemma \ref{lem:IvoryOrth}, the Riemannian metric has the form $a(u,v) du^2 + b(u,v) dv^2$ for some functions $a$ and $b$.
By Lemma \ref{lem:IvoryLem2}, there are functions $\phi(u)$ and $\psi(v)$ such that
\[
\|\phi(u) du + \psi(v) dv\| = 1.
\]
Change the direction of the geodesic $\gamma_+$ and consider a coordinate rectangle with a different ``aspect ratio''.
This yields another pair of functions $\wt\phi(u)$ and $\wt\psi(v)$ with the same property.
Denote $U = \phi^2, V = \psi^2$ and $\wt U = (\wt\phi)^2, \wt V = (\wt\psi)^2$.
Since $\|du\|^2 = \frac1{a}$, etc., the conditions on the norm can be rewritten as
\begin{gather*}
\frac{U}{a} + \frac{V}{b} = 1, \ \ 
\frac{\wt U}{a} + \frac{\wt V}{b} = 1.
\end{gather*}
Solving this system of equations with respect to $a$ and $b$, one obtains

\begin{gather*}
a = \frac{\begin{vmatrix} U & V\\ \wt U & \wt V \end{vmatrix}}{\wt V-V} = 
\left( \frac{U}{\wt U-U} - \frac{V}{\wt V-V} \right) (\wt U-U),\\
b = \frac{\begin{vmatrix} U & V\\ \wt U & \wt V \end{vmatrix}}{U-\wt U} = 
\left( \frac{U}{\wt U-U} - \frac{V}{\wt V-V} \right) (V-\wt V).
\end{gather*}
This shows that the metric is Liouville.
\proofend

\section{Liouville is equivalent to Graves} \label{sect:LtoG}

That Liouville nets have the Graves property, and vice versa, foliations with the Graves property are coordinate lines of a Liouville net, is proved in the Opus Magnum of Darboux  \cite{Dar}, item 589, Livre VI, Chap. I. For a concise exposition of Darboux' arguments we refer the reader to \cite{DraRad}, Section 5.5.
Interestingly, Blaschke's proof of ``Ivory implies Liouville'' is very similar to that of Darboux for ``Graves implies Liouville'': Darboux deduces from the Graves property the existence of two independent distance functions in separated variables.

\appendix
\section{Liouville nets in the plane: Weihnacht-Blaschke theorem}
Here we give a proof of Theorem \ref{thm:LiouvEuc} which uses some of the ideas of \cite{Bla28}, where a $3$-dimensional version of the same result is proved.

Any point $p \in \Omega$ has a neighborhood of the form $(u_{min}, u_{max}) \times (v_{min}, v_{max})$, where $(u, v)$ are any local parameters for the net. It suffices to prove the theorem for such a net rectangle.
We start with the following lemma.

\begin{lemma}
\label{lem:BlaschkeLemma}
For the coordinate lines $\{u = \mathrm{const}\}$ one of the following holds.
\begin{enumerate}
\item
All of them are arcs of conic sections and the map $(u_0, v) \mapsto (u_1, v)$ between any two curves $\{u = u_0\}$ and $\{u = u_1\}$ is the restriction of an affine transformation of $\R^2$.
\item
All of them are straight line segments.
\end{enumerate}
\end{lemma}
\proof
Choose any $u_0, u_1 \in (u_{min}, u_{max})$ and $v_0, v_1, v_2 \in (v_{min}, v_{max})$. For an arbitrary $v \in (v_{min}, v_{max})$ consider the net rectangles $[u_0, u_1] \times [v_i, v]$, $i = 0, 1, 2$. (We abuse the notation, the endpoints of the coordinate intervals might be in a different order.)
Since Liouville nets satisfy the Ivory property, the diagonals in each of these rectangles have equal length:
\begin{multline}
(x(u_1,v) - x(u_0,v_i))^2 + (y(u_1,v) - y(u_0,v_i))^2 \label{eqn:Ivory3Eq}\\
= (x(u_0,v) - x(u_1,v_i))^2 + (y(u_0,v) - y(u_1,v_i))^2, \quad i = 0, 1, 2.
\end{multline}
Here $(x(u,v), y(u,v))$ are Cartesian coordinates of the point with Liouville coordinates $(u,v)$.
Subtract from the equation for $i=0$ the equation for $i = 1$ or $2$. This leads to two linear equations between the Cartesian coordinates of the points $(u_0, v)$ and $(u_1, v)$:
\begin{equation}
\label{eqn:System}
a_i x(u_1,v) + b_i y(u_1,v) + c_i x(u_0,v) + d_i y(u_0,v) + e_i = 0, \quad i = 1, 2.
\end{equation}
The further case distinction depends on the degeneracy of the matrix
\[
\begin{pmatrix}
a_1 & b_1\\
a_2 & b_2
\end{pmatrix}
=
\begin{pmatrix}
x(u_0,v_1) - x(u_0,v_0) & y(u_0,v_1) - y(u_0,v_0)\\
x(u_0,v_2) - x(u_0,v_0) & y(u_0,v_2) - y(u_0,v_0)
\end{pmatrix}.
\]

Since the choices of $u_0, u_1, v_0, v_1, v_2$ are mutually independent, for any choice of $u_0$ one of the following holds.

{\bf Case 1.}
There are $v_0, v_1, v_2$ such that $\det\begin{pmatrix} a_1 & b_1\\ a_2 & b_2 \end{pmatrix} \ne 0$.

Then the linear system \eqref{eqn:System} can be solved for $x(u_1, v)$ and $y(u_1,v)$ which become linear (non-homogeneous) functions of $x(u_0,v)$ and $y(u_0,v)$. Substituting these functions into any of the equations \eqref{eqn:Ivory3Eq} one obtains a quadratic equation on $(x(u_0,v), y(u_0,v))$. It follows that the curve $\{u = u_0\}$ is a conic section or a line, and that for all $u_1$ the map $(u_0,v) \mapsto (u_1,v)$ is the restriction of an affine transformation.

{\bf Case 2.}
For all $v_0, v_1, v_2$ one has $\det\begin{pmatrix} a_1 & b_1\\ a_2 & b_2 \end{pmatrix} = 0$.

For $v_0$ and $v_1$ fixed, this is a linear equation on $x(u_0, v_2)$ and $y(u_0, v_2)$. Thus $\{u = u_0\}$ is a line segment.

If for all $u_0$ Case 2 takes place, then all coordinate curves $\{u = \mathrm{const}\}$ are straight line segments: the second possibility of Lemma \ref{lem:BlaschkeLemma} is realized. If for at least one choice of $u_0$ Case 1 takes place, then all these curves are straight line segments or arcs of conics and are affine images of each other: the first possibility of Lemma \ref{lem:BlaschkeLemma} is realized.
\proofend

If both $u$-coordinate lines and $v$-coordinate lines are straight line segments, then they form an orthogonal grid of lines.
Assume that the curves $\{u = \mathrm{const}\}$ are arcs of conics. From Lemma \ref{lem:BlaschkeLemma} and the orthogonality of the Liouville net we will derive that this family of conics is one of the following.
\begin{itemize}
\item
concentric circles;
\item
confocal ellipses;
\item
confocal hyperbolas;
\item
confocal parabolas with parallel directrices.
\end{itemize}

Denote $\gamma_t = \{u = u_0 + t\}$, and let $\phi_t$ be the affine transformation whose restriction sends $\gamma_0$ to $\gamma_t$ via $(u_0, v) \mapsto (u_0 + t, v)$. 
Take any Cartesian coordinate system $(x,y)$ and write the curve $\gamma_t$ as
\[
\gamma_t = \{(x,y) \mid p^\top Q_t p = 0\},
\]
where
\[
p = \begin{pmatrix} x\\ y\\ 1 \end{pmatrix}, \quad
Q_t = \begin{pmatrix} A_t & b_t\\ b_t^\top & c_t \end{pmatrix},
\]
with $A_t$ a symmertic $2 \times 2$ matrix, $b_t \in \R^2$, $c_t \in \R$.
The affine transformation $\phi_t$ can be written as
\[
\begin{pmatrix} \phi_t(x,y)\\ 1 \end{pmatrix} = \begin{pmatrix} \Phi_t & \psi_t\\ 0 & 1 \end{pmatrix} \begin{pmatrix} x\\ y\\ 1 \end{pmatrix} = R_t p.
\]
Due to $\gamma_t = \phi_t(\gamma_0)$, one has
\[
Q_t = (R_t^{-1})^\top Q_0 R_t^{-1} \Rightarrow \dot Q = \dot S^\top Q_0 + Q_0 \dot S,
\]
where $\dot Q$ and $\dot S$ are the derivatives of $Q_t$ and $R_t^{-1}$ at $t=0$.
Take any point $(x,y) \in \gamma_0$. Then $\gamma_0$ and $t \mapsto \Phi_t \begin{pmatrix} x\\ y \end{pmatrix} + \psi_t$ are the Liouville coordinate curves through $(x,y)$, hence they are orthogonal. One obtains a normal vector to the former by computing the gradient of the function $F(x,y) = p^\top Q_0 p$, and a tangent vector to the latter by taking derivative in $t$. These vectors are collinear:
\[
\dot\Phi \begin{pmatrix} x\\ y \end{pmatrix} + \dot\psi = \alpha(x,y) \left( A_0 \begin{pmatrix} x\\ y \end{pmatrix} + b_0 \right).
\]
(We denote by $\dot\Phi$ and $\dot\psi$ the derivatives of $\Phi_t$ and $\psi_t$ at $t=0$.) The above equation holds for all $(x,y) \in \gamma_0$ for some function $\alpha$ defined on $\gamma_0$. It follows that $\alpha(x,y)$ is constant and
\[
\dot\Phi = \alpha A_0, \quad \dot\psi = \alpha b_0 \quad \Rightarrow \quad \dot R = \alpha
\begin{pmatrix}
A_0 & b_0\\
0 & 0
\end{pmatrix}, \quad
\dot S = -\alpha
\begin{pmatrix}
A_0 & b_0\\
0 & 0
\end{pmatrix},
\]
since $R_0=Id$. 

Assume that the conic $\gamma_0$ is central. In an appropriate Cartesian coordinate system one has
\[
Q_0 =
\begin{pmatrix}
\frac{1}{a} & 0 & 0\\
0 & \frac{1}{b} & 0\\
0 & 0 & 1
\end{pmatrix}
\]
(we do not specify the signs of $a$ and $b$).
This implies
\[
\dot S = -\alpha
\begin{pmatrix}
\frac{1}{a} & 0 & 0\\
0 & \frac{1}{b} & 0\\
0 & 0 & 0
\end{pmatrix}
\]
and
\[
\dot Q = \dot S^\top Q_0+Q_0\dot S= -2\alpha
\begin{pmatrix}
\frac{1}{a^2} & 0 & 0\\
0 & \frac{1}{b^2} & 0\\
0 & 0 & 0
\end{pmatrix}.
\]
At the same time, for the confocal family
\[
Q_t =
\begin{pmatrix}
\frac{1}{a + t} & 0 & 0\\
0 & \frac{1}{b + t} & 0\\
0 & 0 & 1
\end{pmatrix},
\]
one has
\[
\dot Q =
\begin{pmatrix}
- \frac{1}{a^2} & 0 & 0\\
0 & - \frac{1}{b^2} & 0\\
0 & 0 & 0
\end{pmatrix}.
\]
This means that our curve (family) of conics $\gamma_t$ is tangent at $t=0$ to the 
family of confocal conics. Since the same holds at all times $t$, the family $\gamma_t$ is either a confocal family of ellipses or hyperbolas or a concentric family of circles.

Now assume that $\gamma_0$ is a parabola. Choose the coordinates so that the focus of $\gamma_0$ is at the origin, and the axis coincides with the $y$-axis:
\[
Q_0 =
\begin{pmatrix}
a & 0 & 0\\
0 & 0 & 1\\
0 & 1 & -\frac{1}{a}
\end{pmatrix}.
\]
Then
\[
\dot S = -\alpha
\begin{pmatrix}
a & 0 & 0\\
0 & 0 & 1\\
0 & 0 & 0
\end{pmatrix},
\]
and therefore
\[
\dot Q = -2\alpha
\begin{pmatrix}
a^2 & 0 & 0\\
0 & 0 & 0\\
0 & 0 & 1
\end{pmatrix}.
\]
At the same time, for the family of confocal-coaxial parabolas
\[
Q_t =
\begin{pmatrix}
a+t & 0 & 0\\
0 & 0 & 1\\
0 & 1 & -\frac{1}{a+t}
\end{pmatrix},
\]
one has
\[
\dot Q =
\begin{pmatrix}
1 & 0 & 0\\
0 & 0 & 0\\
0 & 0 & \frac{1}{a^2}
\end{pmatrix}.
\]
Thus our curve of conics $\gamma_t$ is tangent at $t=0$ to the confocal-coaxial family of parabolas.

\bigskip
{\bf Acknowledgements}. We are grateful to D. Burago for a consultation. 
AG was supported by RFBR grant 19-51-50005$\_$\YA\FF$\_$\aaa, II by SNCF grants 200021\_169391 and 200021\_179133, and 
ST was supported by NSF grant DMS-1510055.

\end{document}